\documentclass[preprint,12pt]{elsarticle}

\makeatletter
\def\ps@pprintTitle{%
 \let\@oddhead\@empty
 \let\@evenhead\@empty
 \def\@oddfoot{\centerline{\thepage}}%
 \let\@evenfoot\@oddfoot}
\makeatother

\usepackage{lineno,hyperref}
\modulolinenumbers[5]
\usepackage{amssymb}
\usepackage{amsmath}

\begin{document}

\begin{frontmatter}

\title{A new approach to the Pontryagin maximum principle for nonlinear fractional optimal control problems}
\author[myfirstaddress,mysecondaddress,mythirdaddress]{H.M. Ali\corref{mycorrespondingauthor}}
\cortext[mycorrespondingauthor]{Corresponding author}
\ead{hegagi.ali@fc.up.pt}

\author[mysecondaddress]{F. Lobo Pereira}
\author[mythirdaddress]{S. M. A. Gama}

\address[myfirstaddress]{Department of Mathematics, Faculty of Sciences, Aswan University, Aswan, Egypt}
\address[mysecondaddress]{SYSTEC, and Faculty of Engineering, Porto University, Porto, Portugal}
\address[mythirdaddress]{CMUP, and Department of Mathematics, Faculty of Sciences, Porto University, Porto, Portugal}

\begin{abstract}
In this paper, we discuss a new general formulation of fractional optimal control problems whose performance index is in the fractional integral form and the dynamics are given by a set of fractional differential equations in the Caputo sense. We use a new approach to prove necessary conditions of optimality in the form of Pontryagin maximum principle for fractional nonlinear optimal control problems. Moreover, a new method based on a generalization of the Mittag-Leffler function is used to solving this class of fractional optimal control problems. A simple example is provided to illustrate the effectiveness of our main result.
\end{abstract}

\end{frontmatter}

\linenumbers

\section{Introduction}

Fractional optimal control problems (FOCPs) can be regarded as a generalization of classic optimal control problems for which the dynamics of the control system are described by fractional differential equations (FDEs) and might involve a performance index given by fractional integration operator. The reason to formulate and solve FOCPs relies in the fact that there are a significant number of instances in which FDEs describe the behavior of the control systems of interest more accurately than the more common integer differential equations. This is the case, for instance, in diffusion processes, control processing, signal processing, stochastic systems, etc. \cite{key-Das 2011}.

\medspace

Fractional calculus (FC) is a field of Mathematics that deals with integrals and derivatives whose order may be an arbitrary real or complex number, thus generalizing the integer-order differentiation and integration. It started more than $300$ years ago when the notation for differentiation of non-integer order $1/2$ was discussed between Leibnitz and L'Hospital. Since then, fractional calculus has been developed gradually, being now a very active research area of Mathematical Analysis as attested by the vast number of publications (see \cite{key-C G Koha 1990,key-N  Makris1995,key-Mainardi 1996,key-Y  A Rossikin1997}).
There are several different ways of defining fractional derivatives, and, consequently, different types of FOCPs. However, the ones in the sense the Riemann-Liouville and the Caputo have been more widely used. In most of the works that have been published on FOCPs, the state variable is obtained by the Riemann-Liouville or the Caputo fractional integration of the dynamics, but so far, only integer order integral performance indexes have been considered. It also should be noted that several specific numerical techniques have been developed to solve FOCPs. For more details, see \cite{key-Agrawal 2004,key-Agrawal 2010,key-Jelicic 2009,key-Dilfim 2014}.

\medspace

In this paper, we consider FOCPs for which the performance index is given by an integral of fractional order, and the dynamics are mappings specifying the Caputo fractional derivative of the state variable with respect to time. We use Caputo fractional derivatives because it is the most popular one among physicists and scientists. The reason for this is that fractional derivative of constants are zero.  Moreover, the assumptions that we impose on the data of the problem enables a novel approach to the proof based on a generalization of Taylor's expansions and a fractional mean value theorem.
Another contribution of the paper consists on an analytic method to solve the fractional differential equation of the illustrative example based on a generalization of the Mittag-Leffler function and $\alpha$ exponential function.

\medspace

This paper is organized as follows. In the next Section,
we present a brief review of fractional integrals and fractional
derivatives concept and some basic notions specifically pertinent to this work. In Section~\ref{sec:maximum_principle}, we state, discuss and prove necessary conditions of optimality in the form of a Pontryagin Maximum Principle for nonlinear fractional optimal control problems. In Section~\ref{sec:example}, a simple illustrative example of a FOCP solved by a method based on the Mittag-Leffler function is presented. Finally, in Section~\ref{sec:conclusions} we present some conclusions of this research as well as some open challenges.

\section{Some preliminaries in fractional calculus}

There are several definitions of a fractional derivative. In this section, we present a review of some definitions and preliminary facts which are particularly relevant for the results of this article \cite{key-Podlubny 1999,key-Samko Kilbas 1993,key-Kilbas 2006}.

\medskip

\noindent \textbf{Definition 2.1}. Let $f(\cdot)$ be a locally integrable function in interval $[a,b]$. For $t\in[a,b]$ and $\alpha>0$, the left and right Riemann-Liouville fractional integrals are, respectively, defined by

\begin{eqnarray*}
_{a}I_{t}^{\alpha}f(t)&=&\frac{1}{\Gamma(\alpha)}\int_{a}^{t}(t-\tau)^{\alpha-1}f(\tau)d\tau,
\end{eqnarray*}

and

\begin{eqnarray*}
_{t}I_{b}^{\alpha}f(t)&=&\frac{1}{\Gamma(\alpha)}\int_{t}^{b}(\tau-t)^{\alpha-1}f(\tau)d\tau,
\end{eqnarray*}
where $\Gamma(\cdot)$ is the Euler gamma function.

\medskip

\noindent\textbf{Definition 2.2}. Let $f(\cdot)$ be an absolutely continuous function in the interval $[a,b]$. For $t\in[a,b]$ and $\alpha>0$, the left and right Riemann-Liouville fractional derivatives are, respectively, defined by

\begin{eqnarray*}
&&\hspace{-1 cm} _{a}D_{t}^{\alpha}f(t)=\frac{d^{n}}{dt^{n}}\left(_{a}I_{t}^{n-\alpha}f(t)\right)= \frac{1}{\Gamma(n-\alpha)}\left(\frac{d}{dt}\right)^{n} \int_{a}^{t}(t-\tau)^{n-\alpha-1}f(\tau)d\tau,\\
\mbox{and}&&\\
&&\hspace{-1 cm}  _{t}D_{b}^{\alpha}f(t)=(-\frac{d}{dt})^{n}\left(_{t}I_{b}^{n-\alpha}f(t)\right)=
\frac{1}{\Gamma(n-\alpha)}\left(-\frac{d}{dt}\right)^{n} \hspace{-.1 cm}
\int_{t}^{b}\hspace{-.1 cm} (\tau-t)^{n-\alpha-1}f(\tau)d\tau,
\end{eqnarray*}

\noindent where $n\in\mathbb{N}$ is such that $n-1 < \alpha \leq n,$ and $\Gamma(\cdot) $ is as in Definition 2.1.

\medskip

\noindent \textbf{Definition 2.3}. Let $f(\cdot)$ be an integrable continuous function in the $[a,b]$.
For $t\in[a,b]$ and $\alpha>0$, the left and the right Caputo fractional
derivatives are, respectively, defined by

\begin{eqnarray*}
&& \hspace{-.8 cm}
_{a}^{C}\!D_{t}^{\alpha}f(t)=_{a}I_{t}^{n-\alpha}\frac{d^{n}}{dt^{n}}f(t)=
\frac{1}{\Gamma(n-\alpha)}\int_{a}^{t}(t-\tau)^{n-\alpha-1}f^{(n)}(\tau)d\tau,\\
\mbox{and}&&\\
&& \hspace{-.8 cm} _{t}^{C}\!D_{b}^{\alpha}f(t)=_{t}I_{b}^{n-\alpha}\left(-\frac{d}{dt}\right)^{n}f(t)=
\frac{(-1)^{n}}{\Gamma(n-\alpha)}\int_{t}^{b}(\tau-t)^{n-\alpha-1}f^{(n)}(\tau)d\tau,
\end{eqnarray*}

\noindent where $n\in\mathbb{N}$  is such that $n-1 < \alpha \leq n$.

\medskip

\noindent \textbf{Remark 2.1}. If $\alpha=n\in\mathbb{N}_0$, then the Caputo and Riemann-Liouville fractional derivative coincides the ordinary derivative $\displaystyle \frac{d^{n}f(t)}{dt^{n}}$.

\medskip

\noindent\textbf{Remark 2.2}. The Caputo fractional derivative of a constant is always equal to zero. This is not the case with the Riemann-Liouville fractional derivative.

\bigskip

\noindent\textbf{Theorem 2.1} (see \cite{key-Delfim 2013}). Let $\alpha>0$ and $f(\cdot)$ be a differentiable function in $[a,b]$, then
\begin{eqnarray*}
&& _{a}^{C}\!D_{t}^{\alpha}{}_{a}I_{t}^{\alpha}f(t)=f(t), \quad _{t}^{C}\!D_{b}^{\alpha}{}_{t}I_{b}^{\alpha}f(t)=f(t),\\
&&_{a}D_{t}^{\alpha}{}_{a}I_{t}^{\alpha}f(t)=f(t),\quad _{t}D_{b}^{\alpha}{}_{t}I_{b}^{\alpha}f(t)=f(t),\\
\mbox{and}&&\\
&& _{a}I_{b}^{\alpha}{}_{a}^{C}\!D_{t}^{\alpha}f(t)=f(b)-f(a),\quad _{b}I_{a}^{\alpha}{}_{t}^{C}\!D_{b}^{\alpha}f(t)=f(a)-f(b).
\end{eqnarray*}

\bigskip

\noindent\textbf{Theorem 2.2}. Fractional integration by parts.

Let $0<\alpha<1,$  $f(\cdot)$ be a differentiable function in interval
$[a,b]$ and $g(\cdot)\in L_{1}([a,b]).$ Then the following integration by parts formula holds
\begin{eqnarray*}
&&\int_{a}^{b}g(t){}_{a}^{C}\!D_{t}^{\alpha}f(t)dt=
\int_{a}^{b}f(t){}_{t}D_{b}^{\alpha}g(t)dt+[{}_{t}I_{b}^{1-\alpha}g(t)f(t)]_{a}^{b}\\
\mbox{and}&&\\
&&\int_{a}^{b}g(t){}_{t}^{C}\!D_{b}^{\alpha}f(t)dt=
\int_{a}^{b}f(t){}_{a}D_{t}^{\alpha}g(t)dt-[{}_{a}I_{t}^{1-\alpha}g(t)f(t)]_{a}^{b}.
\end{eqnarray*}

\medskip

Another important auxiliary result to prove our Maximum Principle is the generalization of the Bellman-Gronwall Lemma for fractional differential systems. Here, we will consider the following integral from extracted from~\cite{key-lin Gronwall 2013}.

\medskip

\noindent \textbf{Theorem 2.3}. Generalized Bellman-Gronwall inequality.

Suppose $\alpha>0$, $t\in[0,T)$ and the functions $a(t)$, $b(t)$and $u(t)$ are a non-negative and continuous functions on $0\leq t<T$
with
$$u(t)\leq a(t)+b(t)\int_{0}^{t}(t-s)^{\alpha-1}u(s)ds,$$
where $b(t)$ is a bounded and monotonic increasing function on $[0,T)$, then
$$u(t)\leq a(t)+\int_{0}^{t}\left[\sum_{n=1}^{\infty}\frac{(b(t)\Gamma(\alpha))^{n}}{\Gamma(n\alpha)}(t-s)^{n\alpha-1}a(s)\right]ds, \quad t\in[0,T)$$

\medskip

\noindent \textbf{Theorem 2.4}. Generalized Taylor's formula (cf. \cite{key-D Usero,key-Z Odibat}).

\noindent Let $0<\alpha\leq 1$ , $n\in\mathbb{N}$, $f(\cdot)$ be a continuous function in $[a,b]$, $^{C}\!D_{a}^{k\alpha}f(\cdot)\in C[a,b]$ $\forall k=1,\ldots,n$ and $^{C}\!D_{a}^{(n+1)\alpha}f(\cdot)$ is continuous on $[a,b],$ then $\forall x\in[a,b]$  the generalized Taylor's formula for Caputo fractional derivatives is defined by

$$f(x)=\sum_{k=0}^{n}\frac{(x-a)^{k\alpha}}{\Gamma(k\alpha+1)}\; ^{C}\!D_{a}^{k\alpha}f(a)+ R_{n}(x,a),$$

\noindent where
$$R_{n}(x,a)=^{C}\!D_{a}^{(n+1)\alpha}f(\xi)\frac{(x-a)^{(n+1)\alpha}}{\Gamma((n+1)\alpha+1)},$$

\noindent being, for each $x\in[a,b]$, $a\leq\xi\leq x$, and denoting the Caputo fractional derivative of order $\alpha$ by $^{C}\!D_{a}^{\alpha}$.

\noindent Notice that, if $\alpha=1$, the generalized Taylor's formula reduces to the classical Taylor's formula.

\medskip

\noindent \textbf{Lemma 2.1}. (see \cite{key-Hosseinabadi}) Let $f\in C[a,b]$,
$\alpha > 0$, then there exists some $\xi\in(a,b)$ such that

$$I_{a+}^{\alpha}f(x)=\frac{1}{\Gamma(\alpha)}
\int_{a}^{x}(x-t)^{\alpha-1}f(t)dt=f(\xi)\frac{(x-a)^{\alpha}}{\varGamma(1+\alpha)},$$
\noindent where $\xi$ the fractional intermediate value. Remark that there might exist more than
one $\xi$ satisfying this property.

\medskip

\noindent \textbf{Definition 2.4}. The two-parameter Mittag-Leffler function defined by the power series in the form:
$$E_{\alpha,\beta}(z)=\sum_{n=0}^{\infty}\frac{z^{n}}{\Gamma[n\alpha+\beta]}, $$
\noindent where $\alpha$, and $\beta $ are positive parameters. When $\beta=1$, this function is denoted simply by $E_{\alpha}(\cdot)$. We observe that $E_{0,1}(z)=1/(1-z),$
$E_{1,1}(z)=\exp z,$ $E_{1,2}(z)=(\exp z-1)/z,$ and $E_{1,0}(z)=z \exp z.$

Let $A\in\mathbb{R}^{n\times n}$, then the generalization of the two-parameter Mittag-Leffler function becomes
$$E_{\alpha,\beta}(At^{\alpha})=\sum_{n=0}^{\infty}A^{n}\frac{t^{n\alpha}}{\Gamma[n\alpha+\beta]},$$
\noindent and let us define the $\alpha$ exponential matrix function by using Mittag-Leffler function as follows

\begin{equation} e_{\alpha}(A,t)=t{}^{\alpha-1}E_{\alpha,\alpha}(At{}^{\alpha})=
t{}^{\alpha-1}\sum_{n=0}^{\infty}A^{n}\frac{t{}^{n\alpha}}{\Gamma[(n+1)\alpha]}\label{M-L formula}
\end{equation}
The Mittag-Leffler function has several interesting properties. For details see \cite{key-Hegagi 2013 nonlinear, key-Dorota Delfim 2011,key-C. J. Prajapat 2013,key-Z. Wei 2011}.

\section{The FOCP statement and its Maximum Principle}\label{sec:maximum_principle}
In this section, we discuss the FOCP considered in this article, state the associated necessary conditions of optimality, and present its proof which uses an approach that differs from the ones usually adopted for fractional optimal control problems.

Let us consider the simple general problem as follows
\begin{eqnarray}(\bar P) \mbox{ Minimize} &&  _{t_{0}}I_{t_{f}}^\alpha L(t,\bar x(t),u(t)) \nonumber \\
	\mbox{subject to} && _{t_{0}}^{C}\!D_{t}^{\alpha}\bar x(t)=\bar f(t,\bar x(t),u(t)),\quad [t_0,t_f]\;{\cal L}-\mbox{a.e.}\label{eq:1.1}\\
	&& \bar x(t_{0})=\bar x_{0}\in \mathbb{R}^{n}\label{eq:2.2}\\
	&& u(t)\in{\cal U} \label{eq:3.3}
\end{eqnarray}

\noindent where $ {\cal U} =\{ u: [t_0,t_f]\to \mathbb{R}^{m} : u(t)\in \Omega(t)\} $, $\Omega: [t_0,t_f]\to \mathbb{R}^{m}$ is a given set valued mapping, $L: \mathbb{R}^{n}\to \mathbb{R}$  and $\bar f: [t_0,t_f]\times \mathbb{R}^{n}\times\mathbb{R}^{m}\to\mathbb{R}^{n}$  are given functions defining respectively the running cost (or Lagrangian) functional and the fractional dynamics, $_{t_{0}}I_{t_{f}}^\alpha$ is the Riemann-Liouville fractional integral and $_{t_{0}}^{C}\!D_{t}^{\alpha} x $ is the left Caputo fractional derivative of order $\alpha>0$ of the state variable with respect to time.

It is not hard to see that a simple transformation allows us to
convert the problem~$(\bar P)$ into an equivalent one, simply by using this
assumption $_{t_{0}}^{C}\!D_{t}^{\alpha}y(t)=L(t,\bar x(t),u(t))$,
supplemented by the initial condition
$y(t_{0})=0$. Then, we conclude that problem~$(\bar P)$ is equivalent to the one as follows:
\begin{eqnarray}(P) \mbox{ Minimize} && g(x(t_{f}))\nonumber \\
\mbox{subject to} && _{t_{0}}^{C}\!D_{t}^{\alpha}x(t)=f(t,x(t),u(t)),\quad [t_0,t_f]\;{\cal L}-\mbox{a.e.}\label{eq:1} \\
&& x(t_{0})=x_{0}\in \mathbb{R}^{n}\label{eq:2}\\
&& u(t)\in{\cal U}, \label{eq:3}
\end{eqnarray}
\noindent where now $g(x(t_{f}))=y(t_f)$, the state variable $x = col (y,\bar x)$, i.e., it includes $y$ as a first component with initial value at $0$, and the mapping $f= col (L, \bar f)$, i.e., it has $ L $ as first component.

\vspace{.3 cm}

From now on, we consider this as the basic optimal control problem in normal form. We remark that the above problem statement is the simplest one that can be considered containing all the ingredients required for ``bona fide'' optimal control problem.

Now, we will state the assumptions under which our result will be proved.

\begin{itemize}
\item[(H1)] The function $g$ is $C_1$ in $\mathbb{R}^{n}$, i.e., continuously differentiable in its domain.
\item[(H2)] The function $f$ is $C_1$ and Lipschitz continuous with constant $K_f$ in $x$ for all $(t,u)\in \{(t,\Omega(t)): t\in [t_0,t_f]\}$.
\item[(H3)] The function $f$ is continuous in $(t,u)$, for all $x\in \mathbb{R}^{n}$.
\item[(H4)] The set valued map $\Omega:[t_0,t_f]\to \mathbb{R}^{m}$ is compact valued.
\item[(H5)] The set $ f(t,x,\Omega(t))$ is bounded by a certain
  positive constant $M$ for all $(t,x)\in [t_0,t_f] \times \mathbb{R}^{n}$.
\end{itemize}

\noindent These are, by no means, the weakest hypotheses enabling the proof of the maximum principles for FOCPs. However, these ones are of interest in that it allows the particularly simple proof adopted in this article.

\noindent Consider $$H(t,x,p,u) := p^Tf(t,x,u), $$ with $p\in \mathbb{R}^{n}$, to be the Pontryagin function associated to problem $(P)$.

\bigskip

\noindent \textbf{Theorem 3.1 } Let $(x^{*},u^{*})$ be optimal control process for $(P)$. Then, there exists a function $p:[t_{0},t_{f}]\rightarrow\mathbb{R}^{n}$ satisfying
\begin{itemize}
\item the adjoint equation
\begin{equation}
_{t}D_{t_{f}}^{\alpha}p^{T}(t)=p^{T}(t)D_{x}f(t,x^{*}(t),u^{*}(t)),\label{eq:4}
\end{equation}
\item and the transversality condition
\begin{equation}
p^{T}(t_{f})=\nabla_{x}g(x^{*}(t_{f})),\label{eq:5}
\end{equation}
\end{itemize}
\noindent where the operator  $_{t}D_{t_{f}}^{\alpha}$ is right Riemann-Liouville fractional derivative, and $u^{*}:[t_{0},t_{f}]\rightarrow\mathbb{R}^{m}$ is a control strategy such that $u^{*}(t) $ maximizes $[t_0,t_f]$ ${\cal L}$-a.e. the map
$$ u \to H(t,x^{*}(t),p(t), u),$$
\noindent on $\Omega (t)$.

\medskip

\noindent \textbf{Proof of Theorem 3.1}

\medskip

The first key idea is that any perturbation of the optimal control $u^*$ that affects the final value of the state trajectory may increase the cost. Thus, the proof relies on the comparison between the optimal trajectory $x^{*}$ and trajectories $x $ which are obtained by perturbing the optimal control $u^{*}$.
Let $\tau$ be a Lebesgue point in $(t_{0},t_{f})$, and $\varepsilon>0$ sufficiently small so that $ \tau- \varepsilon \geq t_0$. By Lebesgue point in the fractional context, which define in the next definition.

\medskip

\noindent \textbf{Definition 3.1}.  A Lebesgue point of an integrable function $f: \mathbb{R} \rightarrow \mathbb{R}$ is a point $t_0 \in\mathbb{R}$ satisfying

$$\lim_{\epsilon\rightarrow 0^+} \frac{1}{2\varepsilon}\, {_{t_0-\varepsilon}I_{t_0+\varepsilon}^\alpha} \left|f(t)-f(t_0)\right|\rightarrow 0.$$

\noindent It is well known that the subset of Lebesgue points of an integrable function $f$ forms a full Lebesgue measure subset. 

\vspace{0.3cm}

Now, let us consider the perturbed control strategy $u_{\tau,\varepsilon} $ defined by
\begin{eqnarray}
	u_{\tau,\varepsilon}(t)=\left\{\begin{array}{ll}
		\bar{u} & \mbox{if }\; t\in[\tau-\varepsilon,\tau)\\
		u^{*}(t) & \mbox{if }\; t\in[t_{0},t_{f}]\setminus [\tau-\varepsilon,\tau)
	\end{array}\right. \label{ubar}	
\end{eqnarray}

\noindent where $\bar u\in \Omega(t)$ for all $t\in [\tau-\varepsilon,\tau)$, being $\tau$ a Lebesgue point of the reference optimal control strategy. Note that, there is no loss of generality of the choice of $\tau$ due to the fact that the set Lebesgue points is of full Lebesgue measure.

Let $x_{\tau,\varepsilon}$ be the trajectory associated with $u_{\tau,\varepsilon}$, and with $x_{\tau,\varepsilon}(t_0)=x_0$. Clearly, by definition of optimality of $(x^*,u^*)$,

\begin{eqnarray}\left\{\begin{array}{rcl}
0 & \leq & g(x_{\tau,\varepsilon}(t_{f}))-g(x^{*}(t_{f}))\vspace{.2cm}\\
 & = & \nabla_{x}g(x^{*}(t_{f}))[x_{\tau,\varepsilon}(t_{f})-x^{*}(t_{f})] +o(\varepsilon)\vspace{.2cm}\\
 & = & \nabla_{x}g(x^{*}(t_{f}))\Phi_\alpha(t_{f},\tau)[x_{\tau,\varepsilon}(\tau)-x^{*}(\tau)] +o(\varepsilon),\end{array}\right.\label{main_inequality}
\end{eqnarray}

\noindent where $\nabla_{x}g(\cdot)$ is the gradient of $g(\cdot)$, $o(\varepsilon)$ is some positive number satisfying $\displaystyle \lim_{\varepsilon\to 0} \frac{o(\varepsilon)}{\varepsilon} = 0 $, $\Phi_\alpha(\cdot,\cdot)$ is the state transition matrix for the linear fractional differential system $$_{t_{0}}^{C}\!D_{t}^{\alpha}\xi(t)=D_{x}f(t,x^{*}(t),u^{*}(t))\xi(t),$$
\noindent and $x_{\tau,\varepsilon}:[t_{0},t_{f}]\rightarrow\mathbb{R}^{n}$ is the solution to $_{t_{0}}^{C}\!D_{t}^{\alpha}x_{\tau,\varepsilon}(t)=f(t,x_{\tau,\varepsilon}(t), u_{\tau,\varepsilon}(t))$ with $x_{\tau,\varepsilon}(0)=x_{0}$.

\noindent Observe that $x_{\tau,\varepsilon}(t)=x^*(t) $, for all $t\in[t_0,\tau)$.

For all $t\in[\tau-\varepsilon,\tau)$, it is clear that
\begin{eqnarray*} | x_{\tau,\varepsilon}(t)-x^*(t)|&\leq & _{\tau-\varepsilon}I^\alpha_\tau | f(s, x_{\tau,\varepsilon}(s),\bar u)- f(s,x^*(s),u^*(s))|ds\\
&\leq & _{\tau-\varepsilon}I^\alpha_\tau K_f|x_{\tau,\varepsilon}(s)-x^*(s)|ds +2M\frac{\varepsilon^\alpha}{\Gamma(\alpha+1)}\\ &\leq & \frac{\bar M\varepsilon^\alpha}{\Gamma(\alpha+1)},
\end{eqnarray*}
\noindent where $$ \displaystyle \bar M = 2 M \left(1+K_f\sum_{n=1}^\infty
\frac{\Gamma(\alpha)^{n-1}}{ \Gamma(n\alpha+1)}{\varepsilon^{n\alpha}}\right).$$ It is not difficult to show that this series converges and thus $\bar M$ is some finite positive number. The last inequality was obtained by applying Theorem 2.3.

In order to proceed, we need the following auxiliary result.

\medskip

\noindent \textbf{Lemma 3.1}. Consider the general time interval
$[a,b]$ and define the function  $F(t,x)=f(t,x,\bar u),$ where $\bar u$ is like in (\ref{ubar}). Moreover, consider $\tilde x(\cdot)$ and $y(\cdot)$ to be, respectively, solutions to the following fractional differential systems:
\begin{itemize}
\item $_a^C\! D_t^\alpha \tilde x (t) = F(t,\tilde x(t))$ with $ \tilde x(a)=x_a$, and
\item $_a^C\! D_t^\alpha y (t) = D_x F(t,\tilde x(t))y(t) $ with $ y(a)=\bar y \Gamma(\alpha+1)$.
\end{itemize}
\noindent Then, for all $\nu$ positive and sufficiently small real number, we have that $\tilde x_\nu (\cdot)$ solution to the system
$$ _a^C\! D_t^\alpha \tilde x_\nu (t) = F(t,\tilde x_\nu (t)),\quad \tilde x_\nu (a)\in x_a + \nu^\alpha \bar y + o(\nu^\alpha)B_1^n(0), $$
satisfies on $[a,b]$,
$$ \tilde x_\nu (t) \in \tilde x (t) + \frac{\nu^\alpha}{\Gamma(\alpha+1)} y(t) + o(\nu^\alpha)B_1^n(0).$$
\noindent Here, $B_1^n(0)$ denotes the closed unit ball of $\mathbb{R}^{n}$ centered at $0$,

\medskip

\noindent \textbf{Proof of Lemma 3.1}

\noindent After using Taylor\textquoteright s series of fractional order as defined in Theorem 2.4. We conclude next inequality
$$\left|_a^C\!D_t^{\alpha}\left(\tilde x(t)+\frac{\nu^\alpha}{\Gamma(\alpha+1)}y(t) \right)-F(t, \tilde x (t) + \frac{\nu^\alpha}{\Gamma(\alpha+1)}y(t))\right|\leq o(\nu^\alpha)$$

\noindent and we can be extracted this inequality as the follows
$$\left|_a^C\!D_t^{\alpha}\tilde x(t) + _a^C\! D_t^\alpha y (t) \frac{\nu^\alpha}{\Gamma(\alpha+1)}- F(t, \tilde x (t)) - \frac{\nu^\alpha}{\Gamma(\alpha+1)}D_x F(t,\tilde x(t))y(t))\right|\leq o(\nu^\alpha).$$

\noindent Clearly, the first and the third terms cancel each other in the left hand side of the inequality, and, thus, we have
$$\left|_a^C\!D_t^{\alpha}y (t) \frac{\nu^\alpha}{\Gamma(\alpha+1)} - \frac{\nu^\alpha}{\Gamma(\alpha+1)}D_x F(t,\tilde x(t))y(t))\right|\leq o(\nu^\alpha).$$

\noindent By dividing each side by $\nu^{\alpha}$, where
$\frac{o(\nu^{\alpha})}{\nu^{\alpha}} \rightarrow0,$ when $\nu
\rightarrow 0^+,$ we conclude immediately fractional linearized differential system

$$_a^C\! D_t^\alpha y (t) = D_x F(t,\tilde x(t))y(t), \quad  y(a)=\bar y \Gamma(\alpha+1).  {\hskip 3true cm  \square}$$

\medskip

Since, from above, $\displaystyle |x_{\tau,\varepsilon}(\tau)-x^*(\tau)|\leq \bar M \frac{\varepsilon^\alpha}{\Gamma(\alpha+1)}$  for some finite $\bar M$, we may apply Lemma 3.1 to the case of $t\in[\tau, t_f]$.

\noindent By putting $a=\tau$, $\tilde x =x^*$, $y=\xi$, $\nu=\varepsilon$, $\tilde x= x_{\tau,\varepsilon}$ and $$\bar y = f(\tau, x_{\tau,\varepsilon}(\tau), \bar u)- f(\tau, x^*(\tau), u^*(\tau)),$$ Lemma 3.2 readily yields, for almost all $t\in [\tau, t_f]$,
\begin{equation}
x_{\tau,\varepsilon}(t)\in x^*(t)+ \frac{\varepsilon^\alpha}{\Gamma(\alpha+1)} \xi(t) + o(\varepsilon^\alpha)B_1^n(0),\label{perturb-estimate}
\end{equation}
\noindent where $\xi(\cdot)$ satisfies the fractional linearized differential system

\begin{eqnarray}
&&\left\{\begin{array}{rcl}
_\tau^C\! D_t^\alpha \xi (t)&\hspace{-.1 cm}=&\hspace{-.1 cm} D_x f(t, x^*(t),u^*(t)) \xi(t),\qquad {\cal L}-\mbox{a.e. } t\in[\tau,t_f],\vspace{.2 cm}\\
\xi(\tau) &\hspace{-.1 cm}=&\hspace{-.1 cm} f(\tau, x_{\tau,\varepsilon}(\tau), \bar u)- f(\tau, x^*(\tau), u^*(\tau)).\end{array}\right.\label{aux-linear-system}\end{eqnarray}

\noindent By putting together (\ref{perturb-estimate}), and the chain of inequalities in (\ref{main_inequality}) we can immediately write the inequality

\begin{eqnarray}
0&\leq & \nabla_{x}g(x^{*}(t_{f}))\Phi_\alpha(t_{f},\tau) [x_{\tau,\varepsilon}(\tau)-x^{*}(\tau)] +o(\varepsilon)\nonumber\\
&\leq & \frac{\varepsilon^\alpha}{\Gamma(\alpha+1)} \nabla_{x}g(x^{*}(t_{f}))\Phi_\alpha (t_{f},\tau)  \xi(\tau).\label{last-ineq}
\end{eqnarray}

\noindent By putting $\displaystyle p^T(t_f) = -\nabla_x g(x^*(t_f))$ and $\displaystyle p^T(t) = p^T(t_f) \Phi_\alpha(t_{f},t)$, we conclude immediately that the adjoint variable $p:[t_0,t_f]\to \mathbb{R}^{n}$ satisfies the adjoint equation and the transversatility condition, respectively, (\ref{eq:4}) and (\ref{eq:5}).

\noindent This, together with the definition of $\xi(\tau)$ and the
definition of the Pontryagin function, we conclude, after dividing
both sides of the inequality above by $\displaystyle \frac{\varepsilon^\alpha}{\Gamma(\alpha+1)} $, considering the arbitrariness of $\bar u \in \Omega(t)$ and taking the limit as $\varepsilon \to 0^+$ at time $\tau$, that
$$ H(\tau, x^*(\tau), p(\tau), u^*(\tau)) \geq H(\tau, x^*(\tau), p(\tau), \bar u).$$

\noindent The fact that $\tau$ is an arbitrary Lebesgue point in $ [t_0,t_f] $ implies that the maximum condition of our main result holds, that is, $u^{*}(t) $ maximizes, on $\Omega (t)$, the map $ u \to H(t,x^*(t),p(t), u)$, $[t_0,t_f]$ ${\cal L}$-a.e..

Our main result is proved.

\section{Illustrative example}\label{sec:example}

The Pontryagin maximum principle proved in the previous section is now apply to solve a simple problem of resources management that involves minimizing a certain fractional integral subject to given controlled FDEs.

We consider the following problem
\begin{eqnarray}
\mbox{Minimize} && J(u)\label{cost}\\
\mbox{subject to} && _{0}^{C}\!D_{t}^{\alpha}x(t)=u(t)x(t),\quad t\in[0,T], \label{dyn}\\
&& x(0)=x_0,\label{ini}\\
&& u(t)\in[0,1],\label{contr}
\end{eqnarray}
\noindent where $ J(u)=-_{0}I_{T}^{\alpha}(1-u(t))x(t),$ with $0<
\alpha < 1$ and $T >  \Gamma(\alpha+1)^{\alpha^{-1}}.$ Here
$_{0}I_{T}^{\alpha}$ is fractional integral and $_{0}^{C}\!D_{t}^{\alpha}$ is left Caputo fractional derivative.

\noindent The variable $x$ represents a natural resource that takes positive values (note that $x_0 > 0 $ necessarily) ``grows'' according to the law (\ref{dyn}), where the function $u$, designated by control, represents the fraction of the available resource that is used to promote further growth. The overall goal is to find the control strategy that maximizes the amount of accumulated resource over the time interval $[0,T]$ given by the fractional integral (\ref{cost}).

First, we consider an additional state variable component $y$, satisfying
$$_{0}^{C}\!D_{t}^{\alpha}y(t)=(1-u(t))x(t),\quad y(0)=0,$$
\noindent in order obtain the problem statement in the form considered in our main result, that is,
\begin{eqnarray*}
	\mbox{Minimize} && - y(T) \\
	\mbox{subject to} && _{0}^{C}\!D_{t}^{\alpha}x(t)=u(t)x(t), \hspace{1.2 cm} x(0)=x_{0}, \\
	&&_{0}^{C}\!D_{t}^{\alpha}y(t)=(1-u(t))x(t), \quad y(0)=0,\\
    && u(t)\in [0,1].
\end{eqnarray*}
From Theorem 3.1, the adjoint equation (\ref{eq:4}) and the transversality condition~(\ref{eq:5}) for this problem are
\begin{eqnarray}
_{t}D_{T}^{\alpha}p_{1}(t) &\hspace{-.2 cm}=&\hspace{-.2 cm} [p_{1}u^*(t)+p_{2}(1-u^*(t))],\qquad  p_{1}(T)=0,\label{eq:9}\\
_{t}D_{T}^{\alpha}p_{2}(t) &\hspace{-.2 cm}=&\hspace{-.2 cm} 0, \hspace{4.25 cm} p_{2}(T)=1\label{eq:10}
\end{eqnarray}
\noindent where $ _{t}D_{T}^{\alpha}$ is right Riemann-Liouville fractional derivative of order $\alpha$. Thus, we have that $p_{2}(t)\equiv p_{2}(T)=1$, and equation (\ref{eq:9}) becomes
\begin{equation}
_{t}D_{T}^{\alpha}p_{1}(t)=[(p_{1}(t)-1)u^*(t)+1].\label{eq:11}
\end{equation}
From the maximum condition, we know that $u^{*}(t)$ maximizes, ${\cal L}$-a.e. in $[0,1]$, the mapping
$$v\rightarrow p^{T}(t)f(t,x^{*}(t),y^{*}(t),v)=[p_{1}(t)v+p_{2}(t)(1-v)]x^{*}(t).$$

\noindent Since $p_{2}=1$ and $x^{*}(t)>0$ for all $t\in[0,T] $ (this is to conclude from the fact that $x_0 >0$), the mapping to be maximized can be simplified to $v\rightarrow(p_{1}(t)-1)v$. Thus,
given that the system is time invariant, we have that
\begin{eqnarray*} u^{*}(t) &=& \left\{\begin{array}{ll} 1 & \mbox{ if } p_{1}(t)>1\\
 0 & \mbox{ if } p_{1}(t)<1. \end{array}\right.
\end{eqnarray*}
Since $ p_{1}(T)=0$, and $p_1(\cdot)$ is continuous, $\exists b >0$  s.t. $u^{*}(t)=0$ $\forall t\in[T- b,T]$. Thus, from (\ref{eq:9}) we have $_{t}D_{T}^{\alpha}p_{1}(t)=1 $ and, by backwards integration we obtain
\begin{equation}
p_{1}(t)=\frac{(T-t)^{\alpha}}{\Gamma(\alpha+1)}\label{eq:12}
\end{equation}
Obviously that, for $t^{*}=T-(\Gamma(\alpha+1))^{\frac{1}{\alpha}}$, we obtain $p_{1}(t^{*})=1$.
Now, Let us determine the optimal control for $t < t^{*}$. Since, independently of the control $p_1(\cdot)$ remains monotonically decreasing, we have for $t<t^*$, $ u^{*}(t)=1$, and, thus,
\begin{equation}
_{t}D_{T}^{\alpha}p_{1}(t)=p_{1}(t)\label{eq:13}
\end{equation}
The solution of this linear fractional differential equation~(\ref{eq:13}) is given by $p(t)=p(t^*)\Phi_\alpha(t^*,t)$, where $p(t^*)=1 $ and $\Phi_\alpha(t^*,t)$ is the fractional state transition matrix (in fact, scalar-valued) that can be computed by the Mittag-Leffler function defined in the previous section. By setting $\beta=\alpha $, $ A =[1] $ and by replacing $t$ by $t^*-t = T-\Gamma(\alpha+1)^{\alpha^{-1}}-t $, we conclude that
\begin{eqnarray*} p_1(t)&= & e_{\alpha}(1,t^*-t)\\
&=&(t^*-t)^{\alpha-1}E_{\alpha,\alpha}((t^*-t)^\alpha)\\
&=& (t^*-t)^{\alpha-1}\sum_{k=0}^\infty\frac{(t^*-t)^{k\alpha}}{\Gamma((k+1)\alpha)}.
\end{eqnarray*}
Note that if $\alpha=1$, then we have classical solution $e^{T-t-1}$.

Since we have the optimal control $u^*$, we can easily compute the optimal trajectory which satisfies $x^*(0)=x_0 $, and
\begin{eqnarray*}
_{0}^{C}\!D_{t}^{\alpha}x^*(t)&=&
\left\{\begin{array}{ll} x^*(t) &\mbox{ if } t\in [0, t^*]\vspace{.2 cm}\\
0 &\mbox{ if } t\in [t^*, T]\end{array}\right.
\end{eqnarray*}

\noindent We can compute the optimal trajectory $x^* $ by the generalization Mittag-Leffler function, $\forall t \in [0, t^*],$ $x^*(0)=x_0,$ we conclude that
\begin{eqnarray*}
x^*(t)&= & E_{\alpha}(at^\alpha)\\
&=&x_0 E_{\alpha}(t^\alpha)\\
&=& x_0 \sum_{k=0}^\infty\frac{t^{k\alpha}}{\Gamma(k\alpha+1)}.
\end{eqnarray*}

\noindent Note that if $\alpha=1$, then we have classical solution $x_0e^{t}$.

Now we compute the optimal trajectory $x^* $ in the interval $[t^*,T],$ which $u^*=0,$ $x^*(t)=x^*(T),$  we conclude that
\begin{eqnarray*}
x^*(t)&= & x^*(t^*)\\
&=&x_0 E_{\alpha}((t^*)^\alpha)\\
&=& x_0 \sum_{k=0}^\infty\frac{(T-(\Gamma(\alpha+1))^{\frac{1}{\alpha}})^{k\alpha}}{\Gamma(k\alpha+1)}.
\end{eqnarray*}

\noindent Note that if $\alpha=1$, then we have classical solution $x_0e^{T-1}$.

\section{Conclusion}\label{sec:conclusions}

This article concerns the derivation of necessary conditions of optimality in the form of Pontryagin maximum principle for a nonlinear fractional optimal control problem whose differential equation involves the Caputo derivative of the state variable with respect to time. Under mild assumptions on the data of the problem the proof involved the direct application of variational arguments, thus avoiding the often used argument of converting the optimal control problem into a conventional one and, then, express the optimality conditions for this auxiliary problem back in the fractional derivative context. Another interesting novelty consists in the fact that, unlike in most fractional optimal control problem formulations, we consider the cost functional given by a fractional integral of Riemann-Liouville type.

\medspace

A simple example illustrating the application of our maximum principle was presented. The optimal control strategy was computed analytically being the fractional differential adjoint equation solved by using technique based on a generalization Mittag-Leffler function.

\medspace

A natural  sequel of this article concerns the weakening
of the assumptions on the data of the problem. notably the mere
measurability dependence of the dynamics with respect to time and to the control variables. This will certainly require more sophisticated variational arguments and the use of methods and results of nonsmooth analysis. Another direction of research consists in increasing the structure of the fractional optimal control problem by considering additional state endpoint constraints, and state and/or mixed constraints in its formulation. In this case, additional regularity assumptions will be needed to ensure that the obtained necessary conditions of optimality do not degenerate.

\vspace {3mm}

\noindent{\bf Acknowledgements} H.M. Ali is grateful for the financial
support given by Erasmus Mundus-Deusto University through the grant
Erasmus Mundus Fatima Al Fihri Scholarship Program Lot 1 (EMA2 Lot
1). The work of S. Gama was supported, in part, by FCT
through the CMUP -- Centro de Matem\' atica da Universidade
do Porto.


\section*{References}

\begin{thebibliography}{10}
\bibitem{key-Das 2011}S. Das \newline
Functional  Fractional Calculus
(2nd edition), Springer  (2011). 

\bibitem{key-C G Koha 1990}C. G. Koha and J. M. Kelly \newline
Application
of fractional derivatives to seismic analysis of base-isolated models,
Earthquake engineering and Structural dynamics, 19 (1990), pp. 229--241.

\bibitem{key-N Makris1995} N. Makris, G. F. Dargush and
  M. C. Constantinou \newline
Dynamic analysis of viscoelastic-fluid dampers, J. Engrg. Mech.
ASCE,  121 (1995), pp. 1114--1121.

\bibitem{key-Mainardi 1996} F. Mainardi \newline 
Fractional relaxation-oscillation
and fractional diffusion-wave phenomena', Chaos, Solitons and Fractionals,
 7 (1996), pp. 1461--1477.

\bibitem{key-Y A Rossikin1997} Y. A. Rossikin and M. V. Shikova \newline
Application of fractional calculus to dynamic problems of linear
and nonlinear hereditary mechanics of solids, Appl. Mech. Rev, 
50 (1997), pp. 5--67.

\bibitem{key-Agrawal 2004} O. P. Agrawal \newline 
A General formulation
and solution scheme for fractional optimal control problems, Nonlinear
Dynamics,  38 (2004), pp. 323-- 337.

\bibitem{key-Agrawal 2010} O. P. Agrawal, O. Defterli and D. Baleanu \newline
Fractional optimal control problems with several state and control
variables, Journal of Vibration and Control, 16, No. 13 (2010), pp.
1967--1976.

\bibitem{key-Jelicic 2009} Z. D. Jelicic and N. Petrovacki \newline
Optimality
conditions and a solution scheme for fractional optimal control problems,
Struct. Multidiscip. Optim.  38 , No. 6 (2009), pp. 571--581.

\bibitem{key-Dilfim 2014} S. Pooseh, R. Almeida, D. F. M. Torres \newline
Fractional order optimal control problems with free terminal time,
Journal of industrial and management optimization, 10, No. 2 (2014),
pp. 363--381.

\bibitem{key-Podlubny 1999} I. Podlubny \newline
Fractional Differential
Equations, Mathematics in Sciences and Engineering, 198, Academic
Press, San Diego (1999).

\bibitem{key-Samko Kilbas 1993} S. G. Samko, A. A. Kilbas, O. I.
Marichev \newline
Fractional Integrals and Derivatives: Theory and Applications,
Gordon and Breach Science Publishers S.A., Yverdon (1993).

\bibitem{key-Kilbas 2006} A. A. Kilbas, H. M. Srivastava,
  J. J. Trujillo \newline
Theory and Applications of Fractional Differential Equations,
North\textendash Holland Mathematics Studies, 204. Elsevier Science
B. V., Amsterdam (2006).

\bibitem{key-Delfim 2013} M. J. Lazo, and D. F. Torres \newline
The DuBois-Reymond
fundamental lemma of the fractional calculus of variations and an
Euler-Lagrange equation involving only derivatives of Caputo, Journal
of Optimization Theory and Applications,  156, No. 1 (2013), pp. 56--67.

\bibitem{key-lin Gronwall 2013} S. y. Lin \newline
Generalized Gronwall
inequalities and their applications to fractional differential equations,
Journal of Inequalities and Applications, 549, No. 1 (2013).

\bibitem{key-D Usero} D. Usero \newline
Fractional Taylor series for Caputo
fractional derivatives Construction of numerical schemes, Preprint
http://www.fdi.ucm.es/profesor/lvazquez/calcfrac/docs/paperUsero.pdf
(2007).

\bibitem{key-Z Odibat} Z. Odibat and N. Shawagfeh \newline
Generalized
Taylor's formula, Appl. Math. Comput.  186 (2007), pp. 286--293.

\bibitem{key-Hosseinabadi} A. N. Hosseinabadi, and M. Nategh \newline
On fractional mean value, preprint arXiv:1412.6310 (2014).

\bibitem{key-Hegagi 2013 nonlinear} A. A. M. Arafa, S. Z Rida, A.
A. Mohammadein, and H. M. Ali \newline
Solving nonlinear fractional differential
equation by generalized Mittag-Leffler function method, Commun.
Theor.Phys,  59 (2013), pp. 661--663.

\bibitem{key-Dorota Delfim 2011} D. Mozyrska, and D. F. Torres
  \newline
Modified
optimal energy and initial memory of fractional continuous-time linear
systems, Signal Processing, 91, No. 3 (2011), pp. 379--385.

\bibitem{key-C. J. Prajapat 2013} C. J. Prajapati \newline
Certain properties
of Mittag-Leffler function with argument $x^\alpha,\alpha> 0$, Italian Journal of Pure and Applied
Mathematics,  30 (2013), pp. 411-416.

\bibitem{key-Z. Wei 2011} Z. Wei, and W. Dong \newline
Periodic boundary value problems for Riemann
Liouville sequential fractional dierential equations. Electron. J.
Qual. Theory Dier. Equ 87 (2011), p.p. 1-13.

\end{thebibliography}
\end{document}